\documentclass[12 pt a4 paper twoside]{article}
\usepackage{amsmath,amssymb}
\usepackage{hyperref}
\numberwithin{equation}{section}
\newtheorem{theorem}{Theorem}[section]
\newtheorem{lemma}{Lemma}[section]
\newtheorem{definition}{Definition}[section]
\newtheorem{remark}{Remark}[section]

\numberwithin{equation}{section}

\newcommand{\R}{\mathbb{R}}
\newcommand{\N}{\mathbb{N}}

\allowdisplaybreaks
\begin{document}
\begin{center}
\textbf{\textsl{Existence and continuation of solutions of Hilfer fractional differential equations}}\\
\quad

    \textbf{D. B. Dhaigude and $^{*}$Sandeep P. Bhairat}\\
    Department of Mathematics,\\ Dr. Babasaheb Ambedkar Marathwada University,\\ Aurangabad -- 431 004, (M.S) India.\\
    dnyanraja@gmail.com and $^{*}$sandeeppb7@gmail.com
\end{center}
\begin{abstract}
\noindent In this paper we study the existence and continuation of solution to general fractional differential equation with Hilfer fractional derivative. First we establish new local existence theorems. Then we derive the continuation theorems. With the help of continuation theorems derived in this paper, several global existence results are constructed.
\end{abstract}
{theorem,} global solutions.
\footnote{${^{*}}$ Author for correspondence, email: sandeeppb7@gmail.com.}
\section{Introduction}
\paragraph{}Theory of fractional order systems has gained remarkable significance during last few decades due to its real world applications in ostensibly diverse and wide spread fields of applied mathematics, physics and engineering. The monographs \cite{hr,kst,vl,mf,pi} are devoted to such practical problems in control theory, modeling, relaxations and serve as a foundation of fractional order theory in physics and applied sciences. Recently, Hilfer \cite{hr,rh}, Mainardi \cite{mf} discussed various applications of fractional differential equations in their works. Many complex phenomena in nature can be described more accurately using various fractional operators and are characterized by rapid change in their state.

Nowadays, numerous fractional differential operators are present in literature but Riemann-Liouville (R-L) \cite{kst} and Caputo \cite{q,mf} are universally accepted approaches. R-L operator places less constraints on concerned function but fails to physically sound with practically applicable initial conditions. To avoid such a difficulty, scientists accepted Caputo approach which admits lot of properties from classical calculus. But in \cite{dw}, it has shown that, Caputo derivative also has some defects in applications. Concretely, as in \cite{dw}, one has
\begin{equation*}
\lim_{\delta\to0} {^{C}D_{0^+}^{n-\delta}}x(t)=x^{(n)}(t),\qquad
\lim_{\delta\to0} {^{C}D_{0^+}^{n+\delta}}x(t)=x^{(n)}(t)-x^{(n)}(a),\quad \delta>0.
\end{equation*}
Observe that, if $x^{(n)}(a)\neq0$ near an integer $n$, a very small error of measurement on fractional order may result in totally unlike results, which implies a common case of fractional dynamic system when system starts from non-constant state. Such a problem does not arises  in R-L sense. Additionally, when the Caputo fractional derivative is applied to describe the Nutting's law \cite{df,mf,n},
\begin{equation*}
\sigma(t)=\nu D^{k}\epsilon(t),
\end{equation*}
says a constant strain $\epsilon$ implies stress is independent of time $t$, i.e. $\sigma\equiv0.$ This violates the physical properties of real viscoelastic materials. While in the R-L theory, the constant strain $\epsilon$ does not lead to constant stress \cite{sb}. Analogously, we prefer the Hilfer (generalized Riemann-Liouville) derivative operator which interpolates the both R-L as well as Caputo sense.

In the recent investigations, many researchers studied the existence and uniqueness of solution of nonlinear fractional differential equations, see \cite{bal,cj,dr,df,kmf,lv,jd3,pi,ps,zr} and references therein. Infact, the global existence of solution of fractional differential equation is one of the elementary property, see \cite{kst,pi}. In this paper we mainly focus on developing the theory of existence and uniqueness. First we obtain the local existences followed by continuation theorems to extend the existence of solutions globally. The results obtained in this paper generalizes the existing results \cite{cz,ls} in the literature. Equivalently, the works by C. Kou, H. Zhou, C. P. Li \cite{cz} and C. Li, S. Sarwar \cite{ls} and references therein follow as particular cases of our main results.

The rest of the article is organized as follows: in section 2, we collect all the useful definitions and previously known lemmas which are used in construction of our main results. Section 3 devoted to local existence of solutions followed by  continuation results with global existence theorems in section 4. Concluding remarks are given in the last section.

\section{Prerequisites}
This section is devoted to basic definitions and lemmas from \cite{kst} the theory of fractional calculus which are used in subsequent sections. Let $C_{1-\gamma}[0,T]$ is a complete metric space of all continuous functions mapping $[0,T]$ into $\R$ with the metric $d$ defined by \cite{kmf}
\begin{equation*}
  d(x_1,x_2)={\|x_1-x_2\|}_{C_{1-\gamma}[0,T]}:=\max_{t\in[0,T]}|t^{1-\gamma}[x_1(t)-x_2(t)]|,
\end{equation*}
where
\begin{equation*}
C_{1-\gamma}[0,T]=\{x(t):(0,T]\to\R:{t}^{1-\gamma}x(t)\in C[0,T]\}.
\end{equation*}
\begin{definition} \cite{cj}
Let $\Omega=(0,T]$ and $f:(0,\infty)\to\R$ is a real valued continuous function. The Riemann-Liouville fractional integral of a function $f$ of order $\alpha\in{\R}^{+}$ is denoted as $I_{0^+}^{\alpha}f$ and defined by
\begin{equation}\label{d1}
I_{0^+}^{\alpha}f(t)=\frac{1}{\Gamma(\alpha)}\int_{0}^{t}\frac{f(s)ds}{(t-s)^{1-\alpha}},\quad t>0,
\end{equation}
where $\Gamma(\alpha)$ is the Euler's Gamma function.
\end{definition}
\begin{definition} \cite{kst}
Let $\Omega=(0,T]$ and $f:(0,\infty)\to\R$ is a real valued continuous function. The Riemann-Liouville fractional derivative of function $f$ of order $\alpha\in{\R}_{0}^{+}=[0,+\infty)$ is denoted as $D_{0^+}^{\alpha}f$ and defined by
\begin{equation}\label{d2}
D_{0^+}^{\alpha}f(t)=\frac{1}{\Gamma(n-\alpha)}\frac{d^{n}}{dt^{n}}\int_{0}^{t}\frac{f(s)ds}{(t-s)^{\alpha-n+1}},
\end{equation}
where $n=[\alpha]+1,$ and $[\alpha]$ means the integral part of $\alpha,$ provided the right hand side is pointwise defined on $(0,\infty).$
\end{definition}
\begin{definition} \cite{hr} The Hilfer fractional derivative $D_{0^+}^{\alpha,\beta}$ of function $f\in L^{1}(0,T)$ of order $n-1<\alpha<n$ and type $0\leq\beta\leq1$ is defined by
  \begin{equation}\label{d3}
    D_{0^+}^{\alpha,\beta}f(t)=I_{0^+}^{\beta(n-\alpha)}D^{n}I_{0^+}^{(1-\beta)(n-\alpha)}f(t),
  \end{equation}
where $I_{0^+}^{\alpha}$ and $D_{0^+}^{\alpha}$ are Riemann-Liouville fractional integral and derivative defined by \eqref{d1} and \eqref{d2}, respectively.
\end{definition}
\begin{remark}
The Hilfer fractional derivative interpolates between the R-L and Caputo fractional derivative since
\begin{equation*}
 D_{0^+}^{\alpha,\beta}=\begin{cases}
  DI_{0^+}^{1-\alpha}=D_{0^+}^{\alpha}, \quad\quad \beta=0,\\
  I_{0^+}^{1-\alpha}D= {^C D_{0^+}^{\alpha}}, \,\,\, \quad \beta=1.
\end{cases}
\end{equation*}
\end{remark}
Let $0<\alpha<1,0\leq\beta\leq1.$ For the analysis we consider the initial value problem
\begin{equation}\label{a}\begin{cases}
D_{0^+}^{\alpha,\beta}x(t)&=f(t,x),\quad t\in(0,+\infty),\\
I_{0^+}^{1-\gamma}x(0^+)&=x_0, \quad \gamma=\alpha+\beta-\alpha\beta,
\end{cases}
\end{equation}
and the initial value problem (IVP) for the system of differential equations
\begin{equation}\label{s1}\begin{cases}
D_{0^+}^{\alpha,\beta}x_{1}(t)&=f_{1}(t,x_{1},x_{2},..,x_{n}), \\
D_{0^+}^{\alpha,\beta}x_{2}(t)&=f_{2}(t,x_{1},x_{2},..,x_{n}), \\
&\cdots\\
D_{0^+}^{\alpha,\beta}x_{n}(t)&=f_{n}(t,x_{1},x_{2},..,x_{n}), \\
I_{0^+}^{1-\gamma}x_{i}(0^+)&=x_0, \quad \gamma=\alpha+\beta-\alpha\beta, i=1,2,..,n,
\end{cases}
\end{equation}
where $f:{\R}^{+}\times{\R}\to\R$ in IVP \eqref{a}, $f_{i}:{\R}^{+}\times{\R}^{n}\to\R$ in IVP \eqref{s1} have weak singularities with respect to $t$ and satisfies the Lipschitz conditions $$|f(t,x)-f(t,y)|\leq L|x-y|,\,\, L>0,$$
\begin{equation*}
|f_{k}(t,x_{1},x_{2},..,x_{n})-f_{k}(t,y_{1},y_{2},..,y_{n})|\leq\sum_{k=1}^{n}L_{k}|x_{k}-y_{k}|,\,\,L_{k}>0, k=1,2,..,n,
\end{equation*}
respectively to ensure the existence of unique solutions.

Furthermore, the equivalence of Hilfer fractional IVP and it's equivalent integral equation is established in \cite{kmf} the following lemma.
\begin{lemma}\cite{kmf}
Let $\gamma=\alpha+\beta-\alpha\beta$ where $0<\alpha<1$ and $0\leq\beta\leq1.$ Let $f:(0,T]\times\R\to\R$ such that $f(t,x(t))\in C_{1-\gamma}[0,T]$ for any $x\in C_{1-\gamma}[0,T].$ If $x\in C_{1-\gamma}^{\gamma}[0,T],$ then $x$ satisfies IVP \eqref{a} if and only if $x$ satisfies the Volterra fractional integral equation of the second kind
\begin{equation}\label{b}
x(t)=\frac{x_0}{\Gamma(\gamma)}t^{\gamma-1}+\frac{1}{\Gamma(\alpha)}\int_{0}^{t}(t-s)^{\alpha-1}f(s,x(s))ds, \quad t\in(0,\infty).
\end{equation}
\end{lemma}
\begin{lemma}\cite{cz} If $M$ is a subset of $C_{1-\gamma}[0,T].$ Then $M$ is precompact if and only if the following conditions hold:
\item[(i)] $\{t^{1-\gamma}x(t):x\in M\}$ is uniformly bounded,
\item[(ii)] $\{t^{1-\gamma}x(t):x\in M\}$ is equicontinuous on $[0,T].$
\end{lemma}
\begin{lemma}\cite{cj}
Let $a<b<c, 0\leq\mu<1,$ $x\in C_{\mu}[a,b],$ $y\in C[b,c]$ and $x(b)=y(b).$ Define
\begin{equation*}
z(t)=\begin{cases}
x(t), & \mbox{if }\,\, t\in (a,b], \\
y(t), & \mbox{if} \,\,\, t\in [b,c].
\end{cases}
\end{equation*}
Then $z\in C_{\mu}[a,c].$
\end{lemma}
\begin{lemma}
\textbf{Schauder Fixed Point Theorem:} \cite{gd}
Let $U$ be a closed bounded convex subset of a Banach space $X$ and $T:U\to U$ is completely continuous. Then $T$ has a fixed point in $U.$
\end{lemma}
\begin{lemma} \cite{kst}
If $\alpha>0,0\leq\mu<1.$ If $\mu>\alpha,$ then the fractional integrals $I_{0^+}^{\alpha}$ are bounded from $C_{\mu}[0,T]$ into $C_{\mu-\alpha}[0,T].$ If $\mu\leq\alpha,$ then the fractional integrals $I_{0^+}^{\alpha}$ are bounded from $C_{\mu}[0,T]$ into $C[0,T].$
\end{lemma}
\section{Local existence}
In this section, we obtain the local existence of solutions of IVPs \eqref{a} and \eqref{s1}. For this, let us make the following two hypothesis.

$(H_1).$ Let $f:{\R}^{+}\times\R\to\R$ in IVP \eqref{a} be a continuous function and there exists a constant $0\leq\delta<1$ such that $(Ax)(t)=t^{\delta}f(t,x(t))$ is continuous bounded map from $C_{1-\mu}[0,T]$ into $C[0,T],$ where $T$ is positive constant.

$(H_2).$ Let $f_{i}:{\R}^{+}\times{\R}^{n}\to\R$ in IVP \eqref{s1} be a continuous functions and there exists constants $0\leq{\delta}_{i}<1,$ such that $(A_{i}x_{i})(t)=t^{{\delta}_{i}}f_{i}(t,x_{1},x_{2},..,x_{n})$, $i=1,2,..,n$ are continuous bounded map from $C_{1-\mu}[0,T]$ into $C[0,T],$ where $T$ is positive constant.

\begin{theorem}
Suppose that $(H_1)$ hold. Then IVP \eqref{a} has at least one solution $x\in C_{1-\gamma}[0,h]$ for some $(T\geq)h>0.$
\end{theorem}
\textbf{Proof:} Let
\begin{equation*}
E=\bigg\{x\in C_{1-\gamma}[0,T]:{\big\|x-\frac{x_0}{\Gamma(\gamma)}t^{\gamma-1}\big\|}_{C_{1-\gamma}[0,T]}=\sup_{0\leq t\leq T}\big|t^{1-\gamma}x(t)-\frac{x_0}{\Gamma(\gamma)}\big|\leq b\bigg\},
\end{equation*}
where $b>0$ is a constant. Since the operator $A$ is bounded, there exists a constant $M>0,$ such that
$\sup{\big\{|(Ax)(t)|:t\in [0,T],x\in E\big\}}\leq M.$
\begin{equation*}
\hspace{-1cm}\text{Again let}\qquad D_h=\bigg\{x:x\in C_{1-\gamma}[0,h],\sup_{0\leq t\leq h}\big|t^{1-\gamma}x(t)-\frac{x_0}{\Gamma(\gamma)}\big|\leq b\bigg\},
\end{equation*}
where $h=\min\bigg\{{\big(\frac{b\Gamma(\alpha-\delta+1)}{M\Gamma(1-\delta)}\big)}^{\frac{1}{\alpha-\delta}},T\bigg\}.$ Obviously, $D_h\subset C_{1-\gamma}[0,h]$ is nonempty, closed bounded and convex subset.\

Note that $h\leq T,$ define the operator $B$ as follows:
\begin{equation}\label{c}
  (Bx)(t)=\frac{x_0}{\Gamma(\gamma)}t^{\gamma-1}+\frac{1}{\Gamma(\alpha)}\int_{0}^{t}(t-s)^{\alpha-1}f(s,x(s))ds, \quad t\in[0,h].
\end{equation}
It follows from $(H_1)$ and Lemma 2.5 that we have $B(C_{1-\gamma}[0,h])\subset C_{1-\gamma}[0,h].$\

On the other hand by relation \eqref{c}, for any $x\in C_{1-\gamma}[0,h],$ we have
\begin{align*}
 \bigg|t^{1-\gamma}(Bx)(t)-\frac{x_0}{\Gamma(\gamma)}\bigg|&=\bigg|\frac{t^{1-\gamma}}{\Gamma(\alpha)}\int_{0}^{t}(t-s)^{\alpha-1}s^{-\delta}[s^{\delta}f(s,x(s))]ds\bigg|\\
  & \leq\frac{t^{1-\gamma}}{\Gamma(\alpha)}\int_{0}^{t}(t-s)^{\alpha-1}s^{-\delta}Mds \\
  & \leq Mt^{1-\gamma}  {I_{0^+}^{\alpha}(t^{-\delta})}\\
  & \leq \frac{Mh^{\alpha-\gamma-\delta+1}\Gamma(1-\delta)}{\Gamma(\alpha-\delta+1)}\leq b,
\end{align*}
which means $BD_h\subset D_h.$

Next we show that $B$ is continuous. Let $x_n,x\in D_h, {\|x_n-x\|}_{C_{1-\gamma}[0,h]}\to 0$ as $n\to +\infty.$
In view of continuity of $A,$ we have  ${\|Ax_n-Ax\|}_{[0,h]}\to 0$ as $n\to +\infty.$ Now noting that
\begin{align*}
\bigg|t^{1-\gamma}(Bx_n)(t)-t^{1-\gamma}&(Bx)(t)\bigg|=\bigg|\frac{t^{1-\gamma}}{\Gamma(\alpha)}\int_{0}^{t}(t-s)^{\alpha-1}f(s,x_n(s))ds\\
&\hspace{2cm}-\frac{t^{1-\gamma}}{\Gamma(\alpha)}\int_{0}^{t}(t-s)^{\alpha-1}f(s,x(s))ds\bigg|\\
&\leq\frac{t^{1-\gamma}}{\Gamma(\alpha)}\int_{0}^{t}(t-s)^{\alpha-1}s^{-\delta}[s^{\delta}|f(s,x_n(s))-f(s,x(s))|]ds\\
&\leq\frac{t^{1-\gamma}}{\Gamma(\alpha)}\int_{0}^{t}(t-s)^{\alpha-1}s^{-\delta}|(Ax_n)(s)-(Ax)(s)|ds\\
&\leq {\|(Ax_n)(s)-(Ax)(s)\|}_{[0,h]}\frac{t^{1-\gamma}}{\Gamma(\alpha)}\int_{0}^{t}(t-s)^{\alpha-1}s^{-\delta}ds
\end{align*}
we have
\begin{equation*}
{\|(Bx_n)(t)-(Bx)(t)\|}_{C_{1-\gamma}[0,h]}\leq{\|(Ax_n)(s)-(Ax)(s)\|}_{[0,h]}\frac{\Gamma(1-\delta)h^{\alpha-\gamma-\delta+1}}{\Gamma(\alpha-\delta+1)}.
\end{equation*}
Then ${\|(Bx_n)(t)-(Bx)(t)\|}_{C_{1-\gamma}[0,h]}\to 0$ as $n\to+\infty.$ Thus $B$ is continuous. Furthermore, we shall prove that the operator $BD_h$ is continuous. Let $x\in D_h,$ and $0\leq t_1<t_2\leq h.$ For any $\epsilon>0,$ note that
\begin{equation*}
\frac{t^{1-\gamma}}{\Gamma(\alpha)}\int_{0}^{t}(t-s)^{\alpha-1}s^{-\delta}ds=\frac{\Gamma(1-\delta)}{\Gamma(\alpha-\delta+1)}t^{\alpha-\gamma-\delta+1}\to 0 \,\,\text{as}\,\, t\to {0}^{+}, \quad 0\leq\delta<1,
\end{equation*}
there exists a $(h>)\delta_1>0$ such that, for $t\in [0,\delta_1],$
\begin{equation*}
  \frac{2Mt^{1-\gamma}}{\Gamma(\alpha)}\int_{0}^{t}(t-s)^{\alpha-1}s^{-\delta}ds<\epsilon \qquad \text{holds.}
\end{equation*}
In the case with $t_1,t_2\in [0,\delta_1],$ we have
\begin{equation}\label{d}
\begin{split}
\big|&\frac{{t}_{1}^{1-\gamma}}{\Gamma(\alpha)}\int_{0}^{t_1}(t_1-s)^{\alpha-1}f(s,x(s))ds-\frac{t_{2}^{1-\gamma}}{\Gamma(\alpha)}\int_{0}^{t_2}(t_2-s)^{\alpha-1}f(s,x(s))ds\big|\\
&\leq\frac{Mt_{1}^{1-\gamma}}{\Gamma(\alpha)}\int_{0}^{t_1}(t_1-s)^{\alpha-1}s^{-\delta}ds+\frac{Mt_{2}^{1-\gamma}}{\Gamma(\alpha)}\int_{0}^{t_2}(t_2-s)^{\alpha-1}s^{-\delta}ds<\epsilon.
  \end{split}
\end{equation}
In the case with $t_1,t_2\in [\frac{\delta_1}{2},h],$ we get
\begin{align*}
\big|{t}_{1}^{1-\gamma}&(Bx)(t_1)-{t}_{2}^{1-\gamma}(Bx)(t_2)\big|\\
=&\big|\frac{{t}_{1}^{1-\gamma}}{\Gamma(\alpha)}\int_{0}^{t_1}(t_1-s)^{\alpha-1}f(s,x(s))ds-\frac{t_{2}^{1-\gamma}}{\Gamma(\alpha)}\int_{0}^{t_2}(t_2-s)^{\alpha-1}f(s,x(s))ds\big|\\
=&\big|\frac{1}{\Gamma(\alpha)}\int_{0}^{t_1}[{t}_{1}^{1-\gamma}(t_1-s)^{\alpha-1}-t_{2}^{1-\gamma}(t_2-s)^{\alpha-1}] f(s,x(s))ds\\
&\hspace{1cm}-\frac{t_{2}^{1-\gamma}}{\Gamma(\alpha)}\int_{t_1}^{t_2}t_{2}^{1-\gamma}(t_2-s)^{\alpha-1}f(s,x(s))ds\big|
\end{align*}
We see from the fact that if $0\leq\mu_1<\mu_2\leq h,$ then for $0\leq s<\mu_1,$ we have $\mu_{1}^{1-\gamma}(\mu_1-s)^{\alpha-1}>\mu_{2}^{1-\gamma}(\mu_2-s)^{\alpha-1}$ and we obtain
\begin{align*}
\big|\frac{1}{\Gamma(\alpha)}&\int_{0}^{t_1}[{t}_{1}^{1-\gamma}(t_1-s)^{\alpha-1}-t_{2}^{1-\gamma}(t_2-s)^{\alpha-1}] f(s,x(s))ds\big|\\
&\leq \frac{1}{\Gamma(\alpha)}\int_{0}^{t_1}|[{t}_{1}^{1-\gamma}(t_1-s)^{\alpha-1}-t_{2}^{1-\gamma}(t_2-s)^{\alpha-1}]s^{-\delta}|s^{\delta}f(s,x(s))ds\\
&\leq\frac{M}{\Gamma(\alpha)}\int_{0}^{\frac{\delta_1}{2}}|[{t}_{1}^{1-\gamma}(t_1-s)^{\alpha-1}-t_{2}^{1-\gamma}(t_2-s)^{\alpha-1}]s^{-\delta}|ds\\
&\hspace{.4cm}+{({\frac{\delta_1}{2}})}^{-\delta}\frac{M}{\Gamma(\alpha)}\int_{\frac{\delta_1}{2}}^{t_1}[{t}_{1}^{1-\gamma}(t_1-s)^{\alpha-1}-t_{2}^{1-\gamma}(t_2-s)^{\alpha-1}]ds\\
&\leq\frac{2M{(\frac{\delta_1}{2})}^{1-\gamma}}{\Gamma(\alpha)}\int_{0}^{\frac{\delta_1}{2}}(\frac{\delta_1}{2}-s)^{\alpha-1}s^{-\delta}ds\\
&\hspace{.3cm}+\frac{M{(\frac{\delta_1}{2})}^{-\delta}}{\Gamma(\alpha+1)}[{t_2}^{1-\gamma}(t_2-t_1)^{\alpha}-{t_2}^{1-\gamma}(t_2-\frac{\delta_1}{2})^{\alpha}+{t_1}^{1-\gamma}(t_1-\frac{\delta_1}{2})^{\alpha}]\\
&\leq\epsilon+\frac{M{(\frac{\delta_1}{2})}^{-\delta}}{\Gamma(\alpha+1)}[{h}^{1-\gamma}(t_2-t_1)^{\alpha}+{t_2}^{1-\gamma}(t_2-\frac{\delta_1}{2})^{\alpha}+{t_1}^{1-\gamma}(t_1-\frac{\delta_1}{2})^{\alpha}]
\end{align*}
On the other hand,
\begin{align*}
\big|\frac{{t_2}^{1-\gamma}}{\Gamma(\alpha)}\int_{t_1}^{t_2}(t_2-s)^{\alpha-1}f(s,x(s))ds\big|&\leq\frac{{(\frac{\delta_1}{2})}^{-\delta}M}{\Gamma(\alpha)}\int_{t_1}^{t_2}t_{2}^{1-\gamma}(t_2-s)^{\alpha-1}ds\\
&=\frac{{(\frac{\delta_1}{2})}^{-\delta}M}{\Gamma(\alpha+1)}t_{2}^{1-\gamma}(t_2-t_1)^{\alpha}\\
&\leq \frac{{(\frac{\delta_1}{2})}^{-\delta}M{h}^{1-\gamma}}{\Gamma(\alpha+1)}(t_2-t_1)^{\alpha}.
\end{align*}
Clearly, there exist a $\delta,\frac{\delta_1}{2}>\delta>0$ such that, for $t_1,t_2\in[\frac{\delta_1}{2},h],|t_1-t_2|<\delta$ implies
\begin{equation}\label{e}
|{t_1}^{1-\gamma}(Bx)(t_1)-{t_2}^{1-\gamma}(Bx)(t_2)|<2\epsilon.
\end{equation}
It follows from equations \eqref{d} and \eqref{e} that $\{{t}^{1-\gamma}(Bx)(t):x\in D_h\}$ is \textit{equicontinuous.} Obviously, it is clear that $\{{t}^{1-\gamma}(Bx)(t):x\in D_h\}$ is \textit{uniformly bounded} since $BD_h\subset D_h.$ By Lemma 2.2, $BD_h$ is percompact. Therefore $B$ is completely continuous. By Schauder fixed point theorem and Lemma 2.1, the IVP \eqref{a} has a local solution. The proof is thus complete.

\begin{theorem}
Suppose that $(H_2)$ hold. Then IVP \eqref{s1} has at least one solution $x_{i}\in C_{1-\gamma}[0,T]$ for some $(T\geq0)h>0.$
\end{theorem}
\textbf{Proof:} Let
\begin{equation*}
E_{s}=\bigg\{x_{i}\in C_{1-\gamma}[0,T]:{\|x_{i}-\frac{x_0}{\Gamma(\gamma)}t^{\gamma-1}\|}_{C_{1-\gamma}[0,T]}=\sup_{0\leq t\leq T}|t^{1-\gamma}x_{i}-\frac{x_0}{\Gamma(\gamma)}|\leq b_{i}\bigg\},
\end{equation*}
for $b_{i}>0 (i=1,2,..,n)$ are constants. Since the operators $A_{i}, (i=1,2,..,n)$ are bounded then there exists constants $M_{i}>0, (i=1,2,..,n)$ such that $$\sup{ \big\{ |(A_{i}x_{i})(t)|:t\in[0,T], x_{i}\in E_{s}\big\}}\leq M_{i}, i=1,2,..,n.$$
\begin{equation*}
\hspace{-1cm}\text{Again let}\quad D_{ih}=\bigg\{x_{i}:x_{i}\in C_{1-\gamma}[0,h],\sup_{0\leq t\leq h}|t^{1-\gamma}x_{i}(t)-\frac{x_0}{\Gamma(\gamma)}|\leq b_{i}\bigg\},
\end{equation*}
\begin{equation*}
\text{where}\quad h=\min{ \Bigg\{{\bigg(\frac{b_1\Gamma(\alpha-\delta_1+1)}{M_1\Gamma(1-\delta_1)}\bigg)}^{\frac{1}{\alpha-\delta_1}},\cdots, {\bigg(\frac{b_n\Gamma(\alpha-\delta_n+1)}{M_n\Gamma(1-\delta_n)}\bigg)}^{\frac{1}{\alpha-\delta_n}}, \, T \Bigg\}}
\end{equation*}
$\alpha>{\delta}_{i},\,i=1,2,..,n.$ Clearly, $D_{ih}\subset C_{1-\gamma}[0,h]$ is nonempty, closed bounded and convex subsets. Note that $h\leq T,t\in[0,h]$ define operators $B_{i}$ as follows.
\begin{equation}\label{s2}
\begin{cases}
(B_{1}x_{1})(t)&=x_0+\frac{1}{\Gamma(\alpha)}\int_{0}^{t}(t-s)^{\alpha-1}f_{1}(s,x_{1}(s),x_{2}(s),..,x_{n}(s))ds,\\
(B_{2}x_{2})(t)&=x_0+\frac{1}{\Gamma(\alpha)}\int_{0}^{t}(t-s)^{\alpha-1}f_{2}(s,x_{1}(s),x_{2}(s),..,x_{n}(s))ds,\\
&\cdots\\
(B_{n}x_{n})(t)&=x_0+\frac{1}{\Gamma(\alpha)}\int_{0}^{t}(t-s)^{\alpha-1}f_{n}(s,x_{1}(s),x_{2}(s),..,x_{n}(s))ds.
\end{cases}
\end{equation}
By \eqref{s2}, for $x_{i}\in C_{1-\gamma}[0,h],$ we have
\begin{align*}
|t^{1-\gamma}(B_{1}x_{1})(t)-\frac{x_0}{\Gamma(\gamma)}|&=|\frac{t^{1-\gamma}}{\Gamma(\alpha)}\int_{0}^{t}(t-s)^{\alpha-1}s^{-\delta_1}
[s^{\delta_1}f_{1}(s,x_{1}(s),x_{2},..,x_{n}(s))]ds|\\
&\leq\frac{M_1}{\Gamma(\alpha)}t^{1-\gamma} {_{0}I_{t}^{\alpha}(t^{-\delta_1})}=\frac{M_1\Gamma(1-\delta_1)}{\Gamma(\alpha-\delta_1+1)}t^{\alpha-\delta_1-\gamma+1},\\
|t^{1-\gamma}(B_{1}x_{1})(t)-\frac{x_0}{\Gamma(\gamma)}|&\leq \frac{M_1\Gamma(1-\delta_1)}{\Gamma(\alpha-\delta_1+1)}h^{\alpha-\delta_1-\gamma+1}\leq b_{1},\\
|t^{1-\gamma}(B_{2}x_{2})(t)-\frac{x_0}{\Gamma(\gamma)}|&\leq \frac{M_2\Gamma(1-\delta_2)}{\Gamma(\alpha-\delta_2+1)}h^{\alpha-\delta_2-\gamma+1}\leq b_{2},\\
&\cdots\\
|t^{1-\gamma}(B_{n}x_{n})(t)-\frac{x_0}{\Gamma(\gamma)}|&\leq \frac{M_n\Gamma(1-\delta_n)}{\Gamma(\alpha-\delta_n+1)}h^{\alpha-\delta_n-\gamma+1}\leq b_{n},
\end{align*}
which shows that, $B_{i}D_{ih}\subset D_{ih}, i=1,2,..,n.$

Next we show that operators $B_{i}$ are continuous. Let $x_{m},x_{i}\in D_{ih}, m>n$, $i=1,2,..,n$ such that $\|x_{m}-x_{i}\|\to0$ as $m\to+\infty.$
In view of continuity of operators $A_{i},$ we have ${\|A_{i}x_{m}-A_{i}x_{i}\|}_{[0,h]}\to0$ as $m\to+\infty.$ Now noting that
\begin{align*}
|t^{1-\gamma}(B_{i}x_{m})&(t)-t^{1-\gamma}(B_{i}x_{i})(t)|=\\
&|\frac{t^{1-\gamma}}{\Gamma(\alpha)}\int_{0}^{t}(t-s)^{\alpha-1}f_{i}(s,x_{m}(s))ds-\frac{t^{1-\gamma}}{\Gamma(\alpha)}\int_{0}^{t}(t-s)^{\alpha-1}f_{i}(s,x_{i}(s))ds|\\
&\leq \frac{t^{1-\gamma}}{\Gamma(\alpha)}\int_{0}^{t}(t-s)^{\alpha-1}|f_{i}(s,x_{m}(s))-f_{i}(s,x_{i}(s))|ds\\
&\leq \frac{t^{1-\gamma}}{\Gamma(\alpha)}\int_{0}^{t}(t-s)^{\alpha-1}s^{-\delta_{i}}|A_{i}(x_{m})(s)-A_{i}(x_{i})(s)|ds\\
&\leq \frac{t^{1-\gamma}}{\Gamma(\alpha)}\int_{0}^{t}(t-s)^{\alpha-1}s^{-\delta_{i}}ds{\|A_{i}(x_{m})(s)-A_{i}(x_{i})(s)\|}_{[0,h]}.
\end{align*}
\begin{equation*}
\hspace{-1cm}\text{we have}\qquad{\|(B_{i}x_{m})(s)-(B_{i}x_{i})(s)\|}_{[0,h]}\leq \frac{\Gamma(1-\delta_i)}{\Gamma(\alpha-\delta_i+1)}h^{\alpha-\delta_i-\gamma+1}.
\end{equation*}
Then ${\|(B_{i}x_{m})(s)-(B_{i}x_{i})(s)\|}_{[0,h]}\to0$ as $m\to+\infty.$ Thus $B_{i}$ are continuous. Furthermore, we prove that operators $B_{i}D_{ih}$ are continuous. Let $x_{i}\in D_{ih}$ and $0\leq t_1<t_2\leq h.$ For any $\epsilon>0,$ note that
\begin{equation*}
\frac{t^{1-\gamma}}{\Gamma(\alpha)}\int_{0}^{t}(t-s)^{\alpha-1}s^{-\delta_i}ds=\frac{\Gamma(1-\delta_i)}{\Gamma(\alpha-\delta_i+1)}t^{\alpha-\delta_i-\gamma+1}\to0\quad\text{as}\quad t\to{0}^{+},
\end{equation*}
where $0\leq\delta_i<1.$ There exists $\tilde{{\delta}_{i}}>0$ such that for $t\in[0,h],$
\begin{equation*}
\frac{{2M_{i}}t^{1-\gamma}}{\Gamma(\alpha)}\int_{0}^{t}(t-s)^{\alpha-1}s^{-\delta_i}ds<\epsilon
\end{equation*}
holds. In this case, for $t_1,t_2\in[0,\tilde{{\delta}_{i}}],$ we have
\begin{align}\label{e1}
|&\frac{{t_1}^{1-\gamma}}{\Gamma(\alpha)}\int_{0}^{t_1}(t_1-s)^{\alpha-1}f_{i}(s,x_{i}(s))ds-\frac{{t_2}^{1-\gamma}}{\Gamma(\alpha)}\int_{0}^{t_2}(t_2-s)^{\alpha-1}f_{i}(s,x_{i}(s))ds|\\
&\leq\frac{{M_{i}}{t_{1}}^{1-\gamma}}{\Gamma(\alpha)}\int_{0}^{t_1}(t_1-s)^{\alpha-1}s^{-\delta_i}ds+\frac{{M_{i}}{t_{2}}^{1-\gamma}}{\Gamma(\alpha)}\int_{0}^{t_2}(t_2-s)^{\alpha-1}s^{-\delta_i}ds<\epsilon.\nonumber
\end{align}
In the case for $t_1,t_2\in[\frac{\tilde{{\delta}_{i}}}{2},h],$ we get
\begin{align}\label{e2}
|{t_1}^{1-\gamma}&(B_{i}x_{i})(t_1)-{t_2}^{1-\gamma}(B_{i}x_{i})(t_2)|\nonumber\\
=&|\frac{{t_1}^{1-\gamma}}{\Gamma(\alpha)}\int_{0}^{t_1}(t_1-s)^{\alpha-1}f_{i}(s,x_{i}(s))ds-\frac{{t_2}^{1-\gamma}}{\Gamma(\alpha)}\int_{0}^{t_2}(t_2-s)^{\alpha-1}f_{i}(s,x_{i}(s))ds|\nonumber\\
\leq&|\frac{1}{\Gamma(\alpha)}\int_{0}^{t_1}[{t_1}^{1-\gamma}(t_1-s)^{\alpha-1}-{t_2}^{1-\gamma}(t_2-s)^{\alpha-1}]f_{i}(s,x_{i}(s))ds|\nonumber\\
&\hspace{1cm}+|\frac{{t_2}^{1-\gamma}}{\Gamma(\alpha)}\int_{t_1}^{t_2}(t_2-s)^{\alpha-1}f_{i}(s,x_{i}(s))ds|
\end{align}
We see from the fact that if $0\leq\mu_1<\mu_2\leq h,$ then for $0\leq s<\mu_1,$ we have $\mu_{1}^{1-\gamma}(\mu_1-s)^{\alpha-1}>\mu_{2}^{1-\gamma}(\mu_2-s)^{\alpha-1}$ and we obtain from the first term on right hand side of inequality \eqref{e2} that
\begin{align*}
|\frac{1}{\Gamma(\alpha)}\int_{0}^{t_1}[&{t_1}^{1-\gamma}(t_1-s)^{\alpha-1}-{t_2}^{1-\gamma}(t_2-s)^{\alpha-1}]f_{i}(s,x_{i}(s))ds|\\
\leq& |\frac{M_{i}}{\Gamma(\alpha)}\int_{0}^{t_1}[{t_1}^{1-\gamma}(t_1-s)^{\alpha-1}-{t_2}^{1-\gamma}(t_2-s)^{\alpha-1}]s^{-\delta_i}ds|\\
\leq&|\frac{M_{i}}{\Gamma(\alpha)}\int_{0}^{\frac{\tilde{{\delta}_{i}}}{2}}|[{t_1}^{1-\gamma}(t_1-s)^{\alpha-1}-{t_2}^{1-\gamma}(t_2-s)^{\alpha-1}]s^{-\delta_i}|ds\\
&\hspace{0.3cm}+\frac{{(\frac{\tilde{{\delta}_{i}}}{2})}^{-\delta_i}}{\Gamma(\alpha)}\int_{(\frac{\tilde{{\delta}_{i}}}{2})}^{t}|[{t_1}^{1-\gamma}(t_1-s)^{\alpha-1}-{t_2}^{1-\gamma}(t_2-s)^{\alpha-1}]s^{-\delta_i}|ds\\
\leq&|\frac{2M_{i}{(\frac{\tilde{{\delta}_{i}}}{2})}^{1-\gamma}}{\Gamma(\alpha)}\int_{0}^{(\frac{\tilde{{\delta}_{i}}}{2})}(\frac{\tilde{{\delta}_{i}}}{2}-s)^{\alpha-1}s^{-\delta_i}|ds\\
&\hspace{0.3cm}+\frac{M_i{(\frac{\tilde{{\delta}_{i}}}{2})}^{-\delta_i}}{\Gamma(\alpha+1)}[{t_2}^{1-\gamma}(t_2-t_1)^{\alpha}-{t_2}^{1-\gamma}(t_2-\frac{\tilde{{\delta}_{i}}}{2})^{\alpha}+{t_{1}}^{1-\gamma}(t_2-\frac{\tilde{{\delta}_{i}}}{2})^{\alpha}]\\
\leq&\epsilon+\frac{M_i{(\frac{\tilde{{\delta}_{i}}}{2})}^{-\delta_i}}{\Gamma(\alpha+1)}[{h}^{1-\gamma}(t_2-t_1)^{\alpha}+|{t_2}^{1-\gamma}(t_2-\frac{\tilde{{\delta}_{i}}}{2})^{\alpha}+{t_{1}}^{1-\gamma}(t_2-\frac{\tilde{{\delta}_{i}}}{2})^{\alpha}|]
\end{align*}
On the other hand from second term on right hand side of equation \eqref{e2}
\begin{align*}
|\frac{{t_2}^{1-\gamma}}{\Gamma(\alpha)}\int_{t_1}^{t_2}(t_2-s)^{\alpha-1}f_{i}(s,x_{i}(s))ds|
&\leq\frac{M_i{(\frac{\tilde{{\delta}_{i}}}{2})}^{-\delta_i}}{\Gamma(\alpha)}\int_{t_1}^{t_2}{t_2}^{1-\gamma}(t_2-s)^{\alpha-1}ds\\
&=\frac{M_i{(\frac{\tilde{{\delta}_{i}}}{2})}^{-\delta_i}}{\Gamma(\alpha+1)}{t_2}^{1-\gamma}(t_2-t_1)^{\alpha}.
\end{align*}
From above discussion, there exist a $\lambda,(\frac{\tilde{{\delta}_{i}}}{2}>)\lambda>0$ such that for $t_1,t_2\in[\frac{\tilde{{\delta}_{i}}}{2},h]$ and $|t_1-t_2|<\lambda$ implies
\begin{equation}\label{e3}
|{t_1}^{1-\gamma}(B_{i}x_{i})(t_1)-{t_2}^{1-\gamma}(B_{i}x_{i})(t_2)|<2\epsilon.
\end{equation}
It follows from equations \eqref{e1} and \eqref{e3} that $\{t^{1-\gamma}(B_{i}x_{i})(t):x_{i}\in{D_{ih}}\}$ are \textit{equicontinuous.} It is also clear that $\{t^{1-\gamma}(B_{i}x_{i})(t):x_{i}\in{D_{ih}}\}$ is \textit{uniformly bounded} since $B_{i}D_{ih}\subset{D_{ih}}.$ Therefore $B_{i}D_{ih}$ are precompact. So the operators $B_{i}$ are completely continuous. By Schauder fixed point theorem and Lemma 2.1, the IVP \eqref{s1} has a local solution. This completes the proof.

\begin{remark}
If we let $\beta=0$ in IVP \eqref{a}, above Theorem 3.1 yields the local existence (\cite{cz}, Theorem 3.1) associated with R-L IVP (\cite{cz}, equation (1)).
\end{remark}

\begin{remark}
If we let $\beta=1$ in IVP \eqref{a}, above Theorem 3.1 yields the local existence (\cite{ls}, Theorem 3.1) associated with Caputo IVP (\cite{ls}, equation (1.1)).
\end{remark}

\section{Continuation and global existence} In this section, we concerned with continuation of solution of IVP \eqref{a} and then we obtain the global existence. Initially, we need the following definition.
\begin{definition}\cite{cz}
Let $x(t)$ on $(0,\nu)$ and $\tilde{x}(t)$ on $(0,\tilde{\nu})$ both are solutions of IVP \eqref{a}. If $\nu<\tilde{\nu}$ and $x(t)=\tilde{x}(t)$ for $t\in(0,\nu),$ we say that $\tilde{x}(t)$ is continuation of $x(t)$ or $x(t)$ can be continued to $(0,\tilde{\nu}).$ A solution $x(t)$ is noncontinuable if it has no continuation. The existing interval of noncontinuable solution $x(t)$ is called the maximum existing interval of $x(t).$
\end{definition}
\begin{lemma}
\cite{cz} Let $0<\alpha<1,\nu>0,h>0,0\leq\sigma<1, u\in C_\sigma[0,\frac{\nu}{2}]$ and $v\in [\frac{\nu}{2},h].$ Then
\begin{equation*}
I_1(t)=\int_{0}^{\frac{\nu}{2}}(t-s)^{\alpha-1}u(s)ds,\quad I_2(t)=\int_{\frac{\nu}{2}}^{{\nu}}(t-s)^{\alpha-1}v(s)ds
\end{equation*}
are continuous on $[\nu,\nu+h].$
\end{lemma}
\begin{theorem}[Continuation Theorem I]
Assume that $(H_1)$ hold. Then $x=x(t), t\in (0,\nu)$ is noncontinuable if and only if for some $\tau\in (0,\frac{\nu}{2})$ and any bounded closed subset $D\subset[\tau,+\infty)\times\R,$ there exists a $t^{*}\in[\tau,\nu)$ such that $(t^{*},x(t^{*}))\notin D.$
\end{theorem}
\textbf{Proof:} We prove this theorem by contradiction. If possible, suppose that $x=x(t)$ is continuable. Then there exists a solution $\tilde{x}(t)$ defined on $(0,\tilde{\nu}),$ such that $x(t)=\tilde{x}(t)$ for $t\in (0,\nu),$ which implies $\lim_{t\to\nu^{-}}x(t)=\tilde{x}(\nu).$ Define $x(\nu)=\tilde{x}(\nu).$ Evidently, $K=\{(t,x(t)):t\in[\tau,\nu)\}$ is a compact subset of $[\tau,+\infty)\times\R.$ However, there exists no $t^*\in[\tau,\nu)$ such that $(t^*,x(t^*))\notin K.$ This contradiction implies $x(t)$ is noncontinuable.\

We prove converse in two steps. Suppose that there exists a compact subset $\Omega\subset[\tau,+\infty)\times\R,$ such that $\{(t,x(t)):t\in[\tau,\nu)\}\subset\Omega.$ The compactness of $\Omega$ implies $\nu<+\infty.$ By $(H_1),$ there exists a $K>0$ such that $\sup_{(t,x)\in\Omega}|f(t,x)|\leq K.$\\
\textbf{Step:1.} We now show that the $\lim_{t\to\nu^{-}}x(t)$ exists. Let
\begin{equation}\label{g1}
G(s,t)=|\frac{x_0}{\Gamma(\gamma)}s^{\gamma-1}-\frac{x_0}{\Gamma(\gamma)}t^{\gamma-1}|,\quad (s,t)\in[2\tau,\nu]\times[2\tau,\nu],
\end{equation}
\begin{equation}\label{g2}
J(t)=\int_{0}^{\tau}(t-s)^{\alpha-1}s^{-\delta}ds, \quad t\in[2\tau,\nu].
\end{equation}
Easily we can see that $G(s,t)$ and $J(t)$ are uniformly continuous on $[2\tau,\nu]\times[2\tau,\nu]$ and $[2\tau,\nu],$ respectively.

For all $t_1,t_2\in[2\tau,\nu),t_1<t_2,$ by equation \eqref{g1} we have
\begin{align*}
|x(t_1)-x(t_2)|&=|\frac{x_0}{\Gamma(\gamma)}{t_1}^{\gamma-1}-\frac{x_0}{\Gamma(\gamma)}{t_2}^{\gamma-1}+\frac{1}{\Gamma(\alpha)}\int_{0}^{t_1}(t_1-s)^{\alpha-1}f(s,x(s))ds\\
&\hspace{1cm}-\frac{1}{\Gamma(\alpha)}\int_{0}^{t_2}(t_2-s)^{\alpha-1}f(s,x(s))ds|\\
&\leq G(t_1,t_2)+|\frac{1}{\Gamma(\alpha)}\int_{0}^{t_1}(t_1-s)^{\alpha-1}f(s,x(s))ds\\
&\hspace{1cm}-\frac{1}{\Gamma(\alpha)}\int_{0}^{t_2}(t_2-s)^{\alpha-1}f(s,x(s))ds|\\
&\leq G(t_1,t_2)+|\frac{1}{\Gamma(\alpha)}\int_{0}^{\tau}[(t_1-s)^{\alpha-1}-(t_2-s)^{\alpha-1}]s^{-\delta}[s^{\delta}f(s,x(s))]ds|\\
&\hspace{1cm}+|\frac{1}{\Gamma(\alpha)}\int_{\tau}^{t_1}[(t_1-s)^{\alpha-1}-(t_2-s)^{\alpha-1}]f(s,x(s))ds|\\
&\hspace{1cm}+|\frac{1}{\Gamma(\alpha)}\int_{t_1}^{t_2}(t_2-s)^{\alpha-1}f(s,x(s))ds|\\
&\leq G(t_1,t_2)+|\frac{1}{\Gamma(\alpha)}\int_{0}^{\tau}[(t_1-s)^{\alpha-1}-(t_2-s)^{\alpha-1}]s^{-\delta}(Ax)(s)ds|\\
&\hspace{1cm}+\frac{1}{\Gamma(\alpha)}\int_{\tau}^{t_1}[(t_1-s)^{\alpha-1}-(t_2-s)^{\alpha-1}]|f(s,x(s))|ds\\
&\hspace{1cm}+\frac{1}{\Gamma(\alpha)}\int_{t_1}^{t_2}(t_2-s)^{\alpha-1}|f(s,x(s))|ds\\
&\leq G(t_1,t_2)+\frac{{\|Ax\|}_{[0,\tau]}}{\Gamma(\alpha)}\int_{0}^{\tau}[(t_1-s)^{\alpha-1}-(t_2-s)^{\alpha-1}]s^{-\delta}ds\\
&\hspace{1cm}+\frac{K}{\Gamma(\alpha)}\int_{\tau}^{t_1}[(t_1-s)^{\alpha-1}-(t_2-s)^{\alpha-1}]ds+\frac{K}{\Gamma(\alpha)}\int_{t_1}^{t_2}(t_2-s)^{\alpha-1}ds\\
|x(t_1)-x(t_2)|&\leq G(t_1,t_2)+{\|Ax\|}_{[0,\tau]}|J(t_1)-J(t_2)|\\
&\hspace{1cm}+\frac{K}{\Gamma(\alpha+1)}\bigg[\big[(t_1-s)^{\alpha}-(t_2-s)^{\alpha}\big]{\bigg|}_{\tau}^{t_1}+\big[(t_2-s)^{\alpha}\big]{\bigg|}_{t_1}^{t_2}\bigg]\\
&\leq G(t_1,t_2)+{\|Ax\|}_{[0,\tau]}|J(t_1)-J(t_2)|\\
&\hspace{1cm}+\frac{K}{\Gamma(\alpha+1)}\big[2(t_2-t_1)^{\alpha}+(t_1-\tau)^{\alpha}-(t_2-\tau)^{\alpha}\big].
\end{align*}
From the continuity of $G(s,t)$ and $J(t)$ together with the Cauchy convergence criterion, we obtain $\lim_{t\to\nu^{-}}x(t)=x^{*}.$\\
\textbf{Step:2.} Now we show that $x(t)$ is continuable. Since $\Omega$ is closed subset, we have $(\nu,x^{*})\in\Omega.$ Define $x(\nu)=x^{*}.$ Then $x(t)\in C_{1-\gamma}[0,\nu].$ We define operator
\begin{equation*}
  (Ny)(t)=x_1(t)+\frac{1}{\Gamma(\alpha)}\int_{\nu}^{t}(t-s)^{\alpha-1}f(s,x(s))ds, \qquad t\in[\nu,\nu+1],
\end{equation*}
where $y\in[\nu,\nu+1]$ and
\begin{equation*}
  x_1(t)=\frac{x_0}{\Gamma(\gamma)}t^{1-\gamma}\frac{1}{\Gamma(\alpha)}\int_{0}^{\nu}(t-s)^{\alpha-1}f(s,x(s))ds, \qquad t\in[\nu,\nu+1].
\end{equation*}
In view of Lemma 2.5 and Lemma 4.1, we have $N(C[\nu,\nu+1])\subset C[\nu,\nu+1].$ Let
\begin{equation*}
E_b=\big\{ (t,y):\nu\leq t\leq\nu+1, |y|\leq{\max_{\nu\leq t\leq\nu+1}}|x(t)|+b \big\}, \,\, b>0.
\end{equation*}
In view of continuity of $f$ on $E_b,$ we denote $M={\max}_{(t,y)\in E_b}|f(t,y)|.$ Let
\begin{equation*}
E_h=\big\{ y\in[\nu,\nu+h]:{\max_{t\in[\nu,\nu+h]}}|y(t)-x_1(t)|\leq b, y(\nu)=x_1(\nu) \big\},
\end{equation*}
where $h=\min\bigg\{ 1,\big(\frac{\Gamma(\alpha+1)b}{m} \big)^{\frac{1}{\alpha}}\bigg\}.$\\
We claim that $N$ is completely continuous on $E_h.$ First we show the operator $N$ is continuous. In fact, let $\{y_n\}\subseteq C[\nu,\nu+h],$ ${\|y_n-y\|}_{[\nu,\nu+h]}\to 0$ as $n\to+\infty.$ Then we have
\begin{align*}
 |(Ny_n)(t)-(N&y)(t)|=|\frac{1}{\Gamma(\alpha)}\int_{\nu}^{t}(t-s)^{\alpha-1}[f(s,y_n(s))-f(s,y(s))]ds|\\
&\leq\frac{1}{\Gamma(\alpha)}\int_{\nu}^{t}(t-s)^{\alpha-1}|f(s,y_n(s))-f(s,y(s))|ds \\
&\leq {\|f(s,y_n(s))-f(s,y(s))\|}_{[\nu,\nu+h]}\frac{1}{\Gamma(\alpha)}\int_{\nu}^{t}(t-s)^{\alpha-1}ds\\
&\leq {\|f(s,y_n(s))-f(s,y(s))\|}_{[\nu,\nu+h]}\frac{(t-s)^{\alpha}}{\Gamma(\alpha+1)}{\bigg|}_{s=\nu}^{t}\\
&\leq{\|f(s,y_n(s))-f(s,y(s))\|}_{[\nu,\nu+h]}\frac{h^{\alpha}}{\Gamma(\alpha+1)}.
\end{align*}
By virtue of continuity of $f$ on $E_b,$ we have ${\|Ny_n-Ny\|}_{[\nu,\nu+h]}\to 0$ as $n\to+\infty,$ which implies that $N$ is continuous.\

Secondly, we prove that $NE_h$ is euqicontinuous. For any $y\in E_h,$ for which $(Ny)(\nu)=x_1(\nu)$ and
\begin{align*}
 |(Ny)(t)-x_1(t)|&=|\frac{1}{\Gamma(\alpha)}\int_{\nu}^{t}(t-s)^{\alpha-1}f(s,x(s))ds|\\
&\leq\frac{1}{\Gamma(\alpha)}\int_{\nu}^{t}(t-s)^{\alpha-1}|f(s,x(s))|ds \\
&\leq\frac{M}{\Gamma(\alpha)}\int_{\nu}^{t}(t-s)^{\alpha-1}ds\\
&\leq \frac{M(t-s)^{\alpha}}{\Gamma(\alpha+1)}{\bigg|}_{s=\nu}^{t}=\frac{M(t-\nu)^{\alpha}}{\Gamma(\alpha+1)}\\
&\leq\frac{Mh^{\alpha}}{\Gamma(\alpha+1)}\leq b.
\end{align*}
Thus $NE_h\subset E_h.$\

Set $I(t)=\frac{1}{\Gamma(\alpha)}\int_{0}^{\nu}(t-s)^{\alpha-1}f(s,x(s))ds.$ By Lemma 4.1, $I(t)$ is continuous on $[\nu,\nu+h].$ For every $y\in E_h$, $\nu\leq t_1<t_2\leq\nu+h,$ we have
\begin{align}\label{f}
|(Ny)(t_1)-&(Ny)(t_2)|\leq G(t_1,t_2)+\frac{1}{\Gamma(\alpha)}|\int_{0}^{\nu}[(t_1-s)^{\alpha-1}-(t_2-s)^{\alpha-1}]f(s,x(s))ds|\nonumber\\
&\hspace{1cm}+\frac{1}{\Gamma(\alpha)}|\int_{\nu}^{t_1}[(t_1-s)^{\alpha-1}-(t_2-s)^{\alpha-1}]f(s,x(s))ds|\nonumber\\
&\hspace{1cm}+\frac{1}{\Gamma(\alpha)}|\int_{t_1}^{t_2}(t_2-s)^{\alpha-1}f(s,y(s))ds|\nonumber\\
&\leq G(t_1,t_2)+\frac{1}{\Gamma(\alpha)}|\int_{0}^{\nu}[(t_1-s)^{\alpha-1}-(t_2-s)^{\alpha-1}]f(s,x(s))ds|\nonumber\\
&\hspace{1cm}+\frac{1}{\Gamma(\alpha)}\int_{\nu}^{t_1}[(t_1-s)^{\alpha-1}-(t_2-s)^{\alpha-1}]|f(s,x(s))|ds\nonumber\\
&\hspace{1cm}+\frac{1}{\Gamma(\alpha)}\int_{t_1}^{t_2}(t_2-s)^{\alpha-1}|f(s,y(s))|ds\nonumber\\
&\leq G(t_1,t_2)+|I(t_1)-I(t_2)|\nonumber\\
&\hspace{1cm}+\frac{M}{\alpha\Gamma(\alpha)}\bigg[(t_1-s)^{\alpha}-(t_2-s)^{\alpha}{\bigg|}_{\nu}^{t_1}+(t_2-s)^{\alpha}{\bigg|}_{\nu}^{t_1}\bigg]\nonumber\\
&\leq G(t_1,t_2)+|I(t_1)-I(t_2)|\nonumber\\
&\hspace{.5cm}+\frac{M}{\Gamma(\alpha+1)}\bigg[{\big|}(t_1-\nu)^{\alpha}+(t_2-t_1)^{\alpha}-(t_2-\nu)^{\alpha}+(t_2-t_1)^{\alpha}{\big|}\bigg]\nonumber\\
&\leq G(t_1,t_2)+|I(t_1)-I(t_2)|\nonumber\\
&\hspace{.5cm}+\frac{M}{\Gamma(\alpha+1)}\bigg[2(t_2-t_1)^{\alpha}+(t_1-\nu)^{\alpha}-(t_2-\nu)^{\alpha}\bigg]
\end{align}
In view of uniform continuity of $I(t)$ on $[\nu,\nu+h]$ and relation \eqref{f}, we obtain $\{(Ny)(t):y\in E_h\}$ is equicontinuous.
Therefore $N$ is completely continuous. By Schauder's fixed point theorem, $N$ has a fixed point $\tilde{x}(t)\in E_h.$ i.e.
\begin{align*}
\tilde{x}(t)&=x_1(t)+\frac{1}{\Gamma(\alpha)}\int_{\nu}^{t}(t-s)^{\alpha-1}f(s,\tilde{x}(s))ds\\
&=\frac{x_0}{\Gamma(\gamma)}t^{1-\gamma}+\frac{1}{\Gamma(\alpha)}\int_{0}^{t}(t-s)^{\alpha-1}f(s,\bar{x}(s))ds
\end{align*}
\begin{equation*}
\hspace{-5cm}\text{where}\qquad\qquad\bar{x}(t)=\begin{cases}
             x(t), & \mbox{if }\,\, t\in(0,\nu],\\
             \tilde{x}(t), & \mbox{if}\,\,\,t\in[\nu,\nu+h].
           \end{cases}
\end{equation*}
It follows from Lemma 2.3, that $\bar{x}\in C_{1-\gamma}[0,\nu+h]$ and
\begin{equation*}
  \bar{x}(t)=\frac{x_0}{\Gamma(\gamma)}t^{1-\gamma}+\frac{1}{\Gamma(\alpha)}\int_{0}^{t}(t-s)^{\alpha-1}f(s,\bar{x}(s))ds.
\end{equation*}
Therefore, in view of Lemma 2.5, $\bar{x}(t)$ is a solution of IVP \eqref{a} on $(0,\nu+h].$ This yields contradiction since $x(t)$ is noncontinuable. This completes the proof.

Now we present another continuation theorem which is more convenient for application purpose.
\begin{theorem}[Continuation Theorem 2]
Assume that $(H_1)$ hold. Then $x=x(t), t\in(0,\nu),$ is noncontinuable if and only if
\begin{equation}\label{g}
{\lim}_{t\to\nu^{-}}\sup{|M(t)|}=+\infty,
\end{equation}
where $\qquad M(t)=(t,x(t)),\qquad |M(t)|=(t^2+x^2(t))^{\frac{1}{2}}.$
\end{theorem}
\textbf{Proof:} We prove this theorem by contradiction. If possible, suppose $x=x(t)$ is continuable. Then there exists a solution $\tilde{x}(t)$ of IVP \eqref{a} defined on $(0,\tilde{\nu}),$ $\nu<\tilde{\nu},$ such that $x(t)=\tilde{x}(t)$ for $t\in(0,\nu),$ which implies ${\lim}_{t\to\nu^{-}}x(t)=\tilde{x}(\nu).$\\
Thus $|M(t)|\to|M(\nu)|,$ as $t\to\nu^{-},$ which is a contradiction.\

Conversely, suppose that relation \eqref{g} is not true. Then there exists a sequence $\{t_n\}$ and constant $L>0$ such that
\begin{align}\label{h}
 t_n<t_{n+1}&, \,\,n\in \N,\nonumber\\
 {\lim}_{n\to\infty}t_n&=\nu, \,\, |M(t_n)|\leq L,\\
 i.e. \,\, {t_n}^2+x^2(t_n)&\leq {L}^{2}.\nonumber
\end{align}
Since $\{x(t_n)\}$ is bounded convergent subsequence, without loss of generality, set
\begin{equation}\label{i}
{\lim}_{n\to+\infty}x(t_n)=x^{*}.
\end{equation}
Now we show that, for any given $\epsilon>0,$ there exists $T\in(0,\nu),$ such that $|x(t)-x^{*}|<\epsilon, \,\, t\in(T,\nu),$ we have
\begin{equation}\label{j}
{\lim}_{t\to\nu^{-}}x(t_n)=x^{*}.
\end{equation}
For sufficiently small $\tau>0,$ let $E_1=\big\{ (t,x):t\in [\tau,\nu], |x|\leq{\sup}_{t\in[\tau,\nu)}|x(t)|\big\}.$
Since $f$ is continuous on $E_1,$ denote $M={\max}_{(t,y)\in E_1}|f(t,y)|.$
It follows from equations \eqref{h} and \eqref{i} that there exists $n_0$ such that $t_{n_0}>\tau$ and for $n\geq n_0,$ we have $$|x(t_n)-x^{*}|\leq{\frac{\epsilon}{2}}.$$
If equation \eqref{i} is not true, then for $n\geq n_0,$ there exists $\eta_{k}\in(t_n,\nu)$ such that for $t\in(t_n,\eta_n),$ $|x(t)-x^{*}|<\epsilon$ and $|x(\eta_n)-x^{*}|\geq\epsilon$. Thus
\begin{align*}
\epsilon&\leq |x(\eta_n)-x^{*}|\leq |x(t_n)-x^{*}|+|x(\eta_n)-x(t_n)|\\
&\leq\frac{\epsilon}{2}+|\frac{1}{\Gamma(\alpha)}\int_{0}^{t_n}(t_n-s)^{\alpha-1}f(s,x(s))ds-\frac{1}{\Gamma(\alpha)}\int_{0}^{\eta_n}(\eta_n-s)^{\alpha-1}f(s,x(s))ds|\\
&\leq\frac{\epsilon}{2}+\frac{1}{\Gamma(\alpha)}|\int_{0}^{\tau}[(t_n-s)^{\alpha-1}-(\eta_n-s)^{\alpha-1}]f(s,x(s))ds|\\
&\hspace{0.5cm}+\frac{1}{\Gamma(\alpha)}|\int_{\tau}^{t_n}[(t_n-s)^{\alpha-1}-(\eta_n-s)^{\alpha-1}]f(s,x(s))ds|\\
&\hspace{0.5cm}+\frac{1}{\Gamma(\alpha)}|\int_{t_n}^{\eta_n}(\eta_n-s)^{\alpha-1}f(s,x(s))ds|\\
&\leq\frac{\epsilon}{2}+\frac{{\|Ax\|}_{[0,\tau]}}{\Gamma(\alpha)}|J(t_n)-J(\eta_n)|\\
&\hspace{0.5cm}+\frac{M}{\Gamma(\alpha+1)}\bigg[(t_n-s)^{\alpha}-(\eta_n-s)^{\alpha}{\bigg|}_{\tau}^{t_n}+(\eta_n-s)^{\alpha}{\bigg|}_{t_n}^{\eta_n} \bigg]\\
&\leq\frac{\epsilon}{2}+\frac{{\|Ax\|}_{[0,\tau]}}{\Gamma(\alpha)}|J(t_n)-J(\eta_n)|\\
&\hspace{0.5cm}+\frac{M}{\Gamma(\alpha+1)}\bigg[2(\eta_n-t_k)^{\alpha}+(t_n-\tau)^{\alpha}-(\eta_n-\tau)^{\alpha}\bigg],
\end{align*}
where $J(t)$ is defined by \eqref{g2}. In view of the continuity of $J(t)$ on $[t_{{n}_0},\nu],$ for sufficiently large $n\geq n_0,$ we have
\begin{equation*}
\frac{{\|Ax\|}_{[0,\tau]}}{\Gamma(\alpha)}|J(t_n)-J(\eta_n)|+\frac{M}{\Gamma(\alpha+1)}\bigg[2(\eta_n-t_n)^{\alpha}+(t_n-\tau)^{\alpha}-(\eta_n-\tau)^{\alpha}\bigg]<\frac{\epsilon}{2}
\end{equation*}
implies $\epsilon\leq|x(\eta_n)-x^{*}|<\frac{\epsilon}{2}+\frac{\epsilon}{2}=\epsilon.$ This contradicts and ${\lim}_{t\to{\nu^{-}}}x(t)$ exists.\

By using the similar arguments as in the proof of Theorem 4.1, we can easily find the continuation of $x(t).$ The proof is so complete.
\begin{remark}
For $\beta=0,$ Continuation Theorems 1 and 2 reduces to Continuation Theorems for R-L IVP (\cite{cz}, Theorem 4.1, Theorem 4.2, respectively).
\end{remark}
\begin{remark}
For $\beta=1,$ Continuation Theorems 1 and 2 reduces to Continuation Theorems for Caputo IVP (\cite{ls}, Theorem 4.2, Theorem 4.4, respectively).
\end{remark}

Now we study the global existence of solutions for IVP \eqref{a} based on the results obtained in the earlier sections.

Applying Continuation Theorem 2, we can have the following conclusion about the existence of global solution for IVP \eqref{a}.
\begin{theorem}
Suppose that $(H_1)$ hold. Let $x(t)$ is a solution of IVP \eqref{a} on $(0,\nu).$ If $x(t)$ is bounded on $[\tau,\nu)$ for some $\tau>0,$ then $\nu=+\infty.$\
\end{theorem}
Continuing our discussion, we firstly need the following lemma for the further results in our analysis.
\begin{lemma}\cite{cz}
Let $v:[0,b]\to[0,+\infty)$ be a real function, and $w(.)$ be a nonnegative locally integrable function on $[0,b].$ And there exists $a>0$ and $0<\alpha<1,$ such that
\begin{equation*}
  v(t)\leq w(t)+a\int_{0}^{t}(t-s)^{-\alpha}v(s)ds.
\end{equation*}
Then there exists a constant $k=k(\alpha)$ such that for $t\in[0,b],$ we have
\begin{equation*}
  v(t)\leq w(t)+ka\int_{0}^{t}(t-s)^{-\alpha}w(s)ds.
\end{equation*}
\end{lemma}
\begin{theorem}
Suppose that $(H_1)$ hold and there exist three nonnegative continuous functions $l(t),m(t),p(t):[0,\infty)\to[0,\infty)$ such that $|f(t,x)|\leq l(t)m(|x|)+p(t),$ where $m(r)\leq r$ for $r\geq0.$ Then IVP \eqref{a} has one solution in $C_{1-\gamma}[0,\infty).$
\end{theorem}
\textbf{Proof:} The existence of a local solution $x(t)$ of IVP \eqref{a} can be concluded from Theorem 3.1. By Lemma 2.1, $x(t)$ satisfies the integral equation \eqref{b}. Suppose that $[0,\nu),\, \nu<+\infty,$ is the maximum existing interval of $x(t).$ Then
\begin{align*}
|t^{1-\gamma}x(t)|&=|\frac{x_0}{\Gamma(\gamma)}+\frac{t^{1-\gamma}}{\Gamma(\alpha)}\int_{0}^{t}(t-s)^{\alpha-1}f(s,x(s))ds|\\
&\leq\frac{x_0}{\Gamma(\gamma)}+\frac{\nu^{1-\gamma}}{\Gamma(\alpha)}\int_{0}^{t}(t-s)^{\alpha-1}[l(s)m(s^{1-\gamma}|x(s)|)+p(s)]ds\\
&\leq\frac{x_0}{\Gamma(\gamma)}+\frac{\nu^{1-\gamma}}{\Gamma(\alpha)}\int_{0}^{t}(t-s)^{\alpha-1}l(s)m(s^{1-\gamma}|x(s)|)ds\\
&\hspace{1cm}+\frac{\nu^{1-\gamma}}{\Gamma(\alpha)}\int_{0}^{t}(t-s)^{\alpha-1}p(s)ds\\
&\leq\frac{x_0}{\Gamma(\gamma)}+\frac{\nu^{1-\gamma}}{\Gamma(\alpha)}{\|l\|}_{[0,\nu]}\int_{0}^{t}(t-s)^{\alpha-1}m(s^{1-\gamma}|x|)ds\\
&\hspace{1cm}+\frac{\nu^{1-\gamma}}{\Gamma(\alpha)}\int_{0}^{t}(t-s)^{\alpha-1}p(s)ds\\
\end{align*}
take $v(t)=t^{1-\gamma}|x(t)|$,$w(t)=\frac{x_0}{\Gamma(\gamma)}+\frac{\nu^{1-\gamma}}{\Gamma(\alpha)}\int_{0}^{t}(t-s)^{\alpha-1}p(s)ds$ and $a=\frac{\nu^{1-\gamma}{\|l\|}_{[0,\nu]}}{\Gamma(\alpha)}.$

By Lemma 4.2, we know that $v(t)=t^{1-\gamma}|x(t)|$ is bounded on $[0,\nu).$ Thus for any $\tau\in(0,\nu), \, x(t)$ is bounded on $[\tau,\nu).$ By Theorem 4.3, the IVP \eqref{a} has a solution $x(t)$ on $[0,\infty).$

Following result guarantees the existence and uniqueness of global solution of IVP \eqref{a} on ${\R}^{+}.$
\begin{theorem}
Suppose that $(H_1)$ is satisfied and there exists a nonnegative continuous function $l(t)$ defined on $[0,\infty)$  such that $|f(t,x)-f(t,y)|\leq l(t)|x-y|.$ Then IVP \eqref{a} has a unique solution in $C_{1-\gamma}[0,\infty).$
\end{theorem}
The existence of global solution can be obtained by an arguments similar as above. From the Lipschitz-type condition and Lemma 4.2, we can conclude the uniqueness of global solution. We omit the proof here.

\section{Concluding remarks} In this paper, the global existence of a unique solution of nonlinear IVP with Hilfer fractional derivative is proved with the help of fixed point technique and continuation theorems. Continuation theorem 2 is conveniently more applicable in practical problems. Our results in this paper generalizes the existing results in the literature.

\end{document}